\RequirePackage{fix-cm}
\documentclass[smallextended]{svjour3}       
\smartqed  
\usepackage{graphicx}
\usepackage{epsfig,epsfig}
\usepackage{amsmath}
\usepackage{amssymb}
\usepackage{booktabs,arydshln,multirow}
\usepackage{setspace}
\usepackage{xcolor}
\usepackage{stmaryrd}
\usepackage{cite}
\usepackage{graphicx,epstopdf,latexsym,amsfonts,fancyhdr}
 \usepackage[caption=false]{subfig}
 \usepackage[utf8]{inputenc}
\begin{document}

\title{Finite element methods based on two families of second-order numerical formulas for the fractional Cable model with smooth solutions
}

\titlerunning{Novel numerical formulas for the fractional Cable model based on FE method}        

\author{Baoli Yin         \and
        Yang Liu$^*$          \and
        Hong Li           \and
        Zhimin Zhang
}


\institute{* Corresponding author
\\
Baoli Yin \and Yang Liu \and Hong Li \at
              School of Mathematical Sciences, Inner Mongolia University, Hohhot 010021, China. \\
              \email{baolimath@aliyun.com, mathliuyang@imu.edu.cn,  smslh@imu.edu.cn}           
           \and
           Zhimin Zhang \at
              Beijing Computational Science Research Center, Beijing 100193, China.
              \email{zmzhang@csrc.ac.cn}
           \at
           Department of Mathematics, Wayne State University, Detroit, MI 48202, USA.
              \email{zzhang@math.wayne.edu}
}

\date{Received: date / Accepted: date}

\maketitle

\begin{abstract}
We apply two families of novel fractional $\theta$-methods, the FBT-$\theta$ and FBN-$\theta$ methods developed by the authors in previous work, to the fractional Cable model, in which the time direction is approximated by the fractional $\theta$-methods, and the space direction is approximated by the finite element method.
Some positivity properties of the coefficients for both of these methods are derived, which are crucial for the proof of the stability estimates.
We analyse the stability of the scheme and derive an optimal convergence result with $O(\tau^2+h^{r+1})$ for smooth solutions, where $\tau$ is the time mesh size and $h$ is the spatial mesh size.
Some numerical experiments with smooth and nonsmooth solutions are conducted to confirm our theoretical analysis.
To overcome the singularity at initial value, the starting part is added to restore the second-order convergence rate in time.
\end{abstract}

\section{Introduction}
\label{intro}
In recent years, the fractional differential equations (FDEs) have attracted much attention for its faithfully reflecting the phenomenons in science and engineering, such as in biology, physics, control system and ecology, see \cite{Metzler,LijcHYP,Magin,PODLUBNY,Hesthaven,YanFord,ZhangJiang,Hassani}.
Among these FDEs, the fractional Cable equations introduced by Henry and Langlands \cite{Henry1} are to model electrodiffusion of ions in nerve cells.
The finite domain solutions of the linear fractional Cable equation were derived by Langlands et at. \cite{Langlands} with the generalized Mittag-Leffler function.
Nonetheless, efficient numerical algorithms are needed to effectively derive the approximate solutions since the expression of the analytic ones is complicated.
Liu et al. \cite{QQ} proposed two implicit numerical algorithms for numerically solving the fractional Cable problem within the finite difference framework. Liu et al. \cite{Zhengguang} considered the L1 method when discretizing the fractional Cable model in temporal direction, and presented a fast solution technique to accelerate Toeplitz matrix-vector multiplications arising
from finite difference discretization.
In \cite{Zhangzm1,Liuydyw2,Liud1,Wangyj1,Zhuangph1,Maskari}, some authors developed the finite element method considering different ideas for the fractional Cable model and gave the detailed numerical analysis on convergence.
Lin et al. \cite{Xucj6} developed spectral methods for the fractional Cable model.
Yang et al. \cite{Jiang} applied the time-space spectral Legendre tau method to the direct problem.
The element free Galerkin technique was also developed by Dehghan and Abbaszadeh \cite{Dehghan} for the fractional Cable model with a Dirichlet boundary condition.
Zheng and Zhao \cite{Zhengang} analyzed the fractional Cable equation by the discontinuous Galerkin finite element method.
\par
From the methods above one can see that the key point of efficiently deriving the numerical solutions is developing efficient methods to discretize the fractional derivative of the equation, and theoretically showing that the resulted scheme is stable with a high-order convergence rate.
To this end, some high-order approximation formulas were developed for the fractional calculus, see \cite{Lubich1,Sun1,Alikhanov,McLean1,Ding,Jinbt1,Banjai,Liuzll,Duywly,LiuYin2,YanFord,FengZhuang}.
As is well known that the solutions of fractional PDEs show some singularity at the initial value \cite{Martin2}, some methods or techniques were developed to cope with such difficulty, see \cite{Lubich1,Martin,Zeng1,Jin2}.
In this paper, we apply two families of novel fractional $\theta$-methods, i.e., the fractional BT-$\theta$ (FBT-$\theta$) method and fractional BN-$\theta$ (FBN-$\theta$) method (see the generating functions for both of these two methods defined by (\ref{N.2})), developed by authors in \cite{LiuYin}, to the fractional Cable model,
\begin{equation}\label{I.1}
\begin{cases}
  u_t={}_{RL}D_{0,t}^{1-\gamma}\Delta u-\mu^2{}_{RL}D_{0,t}^{1-\kappa}u+f(\boldsymbol z,t), & (\boldsymbol{z},t)\in\Omega\times (0,T], \\
  u(\boldsymbol{z},t)=0,  &(\boldsymbol{z},t)\in\partial{\Omega}\times [0,T], \\
  u(\boldsymbol{z},0)=u_0(\boldsymbol{z}), & \boldsymbol{z}\in \bar{\Omega}=\Omega \cup \partial{\Omega},
\end{cases}
\end{equation}
where $T>0$, $\gamma$, $\kappa \in (0,1)$, $\Omega \subset \mathbb{R}^d$ is bounded spatial interval (when $d=1$) or convex polygonal spatial domain (when $d=2$), respectively.
$f$ and $u_0$ are given smooth functions.
${}_{RL}D_{0,t}^{\alpha}$ denotes the Riemann-Liouville fractional derivative in time of order $\alpha$ defined by
\begin{equation}\label{I.2}
{}_{RL}D_{0,t}^{\alpha}\psi=\frac{1}{\Gamma(1-\alpha)}\frac{\partial}{\partial t}\int_{0}^{t}\frac{\psi(s)}{(t-s)^{\alpha}}\mathrm{d}s,\quad \alpha \in (0,1),
\end{equation}
where $\Gamma(z)$ denotes the Gamma function.
\par
Our contributions in this paper mainly focus on two aspects:
\par
$\bullet$ Some positivity properties (\ref{S.2}) of the coefficients of the fractional $\theta$-methods are derived which are crucial for the analysis of the stability of the numerical scheme.
Optimal error estimates are derived for smooth solutions. Further, we develop the estimate (\ref{S.3}) with the tool of generating functions.
\par
$\bullet$ Solutions with weak regularity are tested for the model (\ref{I.1}) when applying the fractional $\theta$-methods with a starting part.
The optimal second-order convergence rate in time is obtained.
\par
The outline of the rest of the paper is as follows:
In section 2, we state the novel fractional-$\theta$ methods from the aspect of generating functions and give some recursive formulas to efficiently get the convolution weights. Based on the finite element method in space direction, the fully discrete scheme of (\ref{I.1}) is derived.
In section 3, we first prove some positivity properties of the coefficients of the fractional-$\theta$ methods and then derive the stable estimates of our schemes.
Section 4 mainly focus on the analysis of the error estimates, and the optimal convergence result $O(\tau^2+h^{r+1})$ is obtained for smooth solutions.
In section 5, we implement some numerical experiments to further confirm our theoretical analysis. For the one-dimensional example solutions with weak regularity are tested with some correction terms added.
For the two-dimensional example, we assume the solution is sufficiently smooth such that only the convolution part is needed to approximate the derivatives in the equation.
Finally, we make some conclusions in section 6 and discuss some techniques may be useful for the fractional $\theta$-methods when applied to other types of fractional PDEs.
\par
Throughout the article, we denote by $\|\cdot\|$ the norm in $L^2(\Omega)$ space, and define $\|\cdot\|_m$ with $m\in \mathbb{Z}$ as the $H^2(\Omega)$ space norm. Hence we have $\|\cdot\|_0=\|\cdot\|$ by the definition of $\|\cdot\|_m$.
The generic constants $C>0$ may be different at different occurrence, independent of time mesh $\tau$ and spatial mesh $h$.
\section{Numerical schemes}
To derive the numerical schemes of the fractional Cable model (\ref{I.1}), we first divide the temporal interval $[0,T]$ into equally separated intervals:
$0=t_0<t_1<\cdots<t_n<\cdots<t_N=T$ with $t_n=n\tau$ where $\tau:=T/N$.
For brevity, denote $u(t_n)$ by $u^n$.
For a sequence $\{\omega_k\}_{k=0}^{\infty}$ we identify it with its generating power series $\omega(\xi)=\sum_{k=0}^{\infty}\omega_k \xi^k$, and viceversa. Under proper conditions, $\omega(\xi)$ actually defines a function of $\xi$, i.e., the generating function $\omega(\xi)$.
\par
Define the discrete fractional operator $D_{\tau,\omega}^{\alpha}$ as:
\begin{equation}\label{N.1}
 D_{\tau,\omega}^{\alpha} \varphi^n:=\tau^{-\alpha}\sum_{k=0}^{n}\omega^{(\alpha)}_{n-k}\varphi^k
 +\tau^{-\alpha}\sum_{j=1}^{s}\omega^{(\alpha)}_{n,j}\varphi^j,
\end{equation}
where the convolution weights $\omega^{(\alpha)}_k$ in the convolution part $\tau^{-\alpha}\sum_{k=0}^{n}\omega^{(\alpha)}_{n-k}\varphi^k$ are the coefficients defined by some generating functions.
The starting weights $\{\omega^{(\alpha)}_{n,j}\}_{j=1}^s$ in the starting part $\tau^{-\alpha}\sum_{j=1}^{s}\omega^{(\alpha)}_{n,j}\varphi^j$ are derived by letting
\begin{equation}\label{N.1.0}
{}_{RL}D^{\alpha}_{0,t} t^{\sigma_i}|_{t=t_n}=D_{\tau,\omega}^{\alpha}t^{\sigma_i}_n
\end{equation}
exactly hold for $i=1,\cdots,s$ (see \cite{Lubich1} and \cite{LiuYin}),  where we have assumed that the solution of (\ref{I.1}) can be expanded at initial time with the expression (see \cite{Langlands})
\begin{equation}\label{N.1.1}
u(t)-u(0)=u^{(1)}(t)+u^{(2)}(t), \quad u^{(1)}=\sum_{l=1}^{r}c_l t^{\sigma_l}, \quad u^{(2)}=c_{r+1}t^{\sigma_{r+1}}+\zeta(t)t^{\sigma_{r+2}},
\end{equation}
and $\sigma_l$'s satisfy ${\color{red}0}<\sigma_1<\sigma_2<\cdots<\sigma_r<\sigma_{r+1}<\sigma_{r+2}$, $\zeta(t)$ is a smooth function over $[0,T]$.
We note that for $\sigma_1 \geq 3$, which means the solution is smooth enough at the origin for our schemes, we can omit the starting part in the approximation formula (\ref{N.1}) (see Example 2 in section 5).
However, for solutions with weak regularity at initial value, the starting part is crucial to recovering a second-order convergence rate (see Example 1 in section 5).
\par
In the following discussions we mainly analyse two families of novel fractional $\theta$-methods applied to the equation (\ref{I.1}), which, from the aspect of generating function, can be stated as (see \cite{LiuYin}),
\begin{equation}\label{N.2}
\begin{split}
\text{FBT-$\theta$ method:}
\\
\omega^{(\alpha)}(\xi)=&(1-\theta+\theta \xi)^{-\alpha}[(3/2-\theta)-(2-2\theta)\xi+(1/2-\theta)\xi^2]^{\alpha},
\\
\text{FBN-$\theta$ method:}
\\
\omega^{(\alpha)}(\xi)=&(1+\alpha\theta-\alpha\theta \xi)\big[(3/2-\theta)-(2-2\theta)\xi+(1/2-\theta)\xi^2\big]^{\alpha},
\end{split}
\end{equation}
with $\theta \in (-\infty,\frac{1}{2})$ and $\theta \in [-\frac{1}{2\alpha},1]$, respectively.
\par
We note that when taking $\alpha=1$ and $\theta=0$, both of the methods reduce to the approximation for the first derivative by the BDF2.
And for $\theta=0$, both of the methods coincide with the fractional BDF2;
for $\theta=\frac{1}{2}$, the FBN-$\theta$ method becomes the generalized Newton-Gregory formula.
Several papers examined the special cases mentioned above, for example, see \cite{Peng,Jin1,Jin2}.
\par
According to the appendix in \cite{LiuYin}, we can obtain the convolution weights $\{\omega_k^{(\alpha)}\}_{k=0}^N$ by a recursive formula whose algorithm complexity is of $O(N)$.
We state the algorithm in the following lemmas.
\begin{lemma}\label{NC.1}(See \cite{LiuYin})
The convolution weights $\omega_k$ which are defined as the coefficients of the generating function for the FBT-$\theta$ method can be derived by the recursive formula
\begin{equation}\label{NC.1.1}\begin{split}
\omega_0=&\bigg(\frac{3-2\theta}{2-2\theta}\bigg)^{\alpha},
\quad
\omega_1=\frac{\phi_0\omega_0}{\psi_0},
\quad
\omega_2=\frac{1}{2\psi_0}[(\phi_0-\psi_1)\omega_1+\phi_1\omega_0],
\\
\omega_{k}=&\frac{1}{k\psi_0}\sum_{j=1}^{3}[\phi_{j-1}-(k-j)\psi_{j}]\omega_{k-j}, \quad k \geq 3,
\end{split}\end{equation}
where,
\begin{equation}\label{NC.1.2}
\phi_0=-\frac{\alpha}{2}(2\theta^2-5\theta+4),
\quad \phi_1=-\alpha(2\theta-1)(1-\theta),
\quad \phi_2=-\frac{\alpha\theta}{2}(2\theta-1),
\end{equation}
and
\begin{equation}\label{NC.1.3}\begin{split}
\psi_0=\frac{1}{2}(3-2\theta)(1-\theta),
~\psi_1=\frac{1}{2}(1-2\theta)(3\theta-4),
\\
\psi_2=\frac{1}{2}(1-\theta)(1-6\theta),
~\psi_3=\frac{1}{2}\theta(1-2\theta).
\end{split}\end{equation}
\end{lemma}
\begin{lemma}\label{NC.2}(See \cite{LiuYin})
The convolution weights $\omega_k$ which are defined as the coefficients of the generating function for the FBN-$\theta$ method can be derived by the recursive formula
\begin{equation}\label{NC.2.1}\begin{split}
\omega_0=&2^{-\alpha}(1+\alpha\theta)(3-2\theta)^{\alpha},
\quad
\omega_1=\frac{\phi_0\omega_0}{\psi_0},
\quad
\omega_2=\frac{1}{2\psi_0}[(\phi_0-\psi_1)\omega_1+\phi_1\omega_0],
\\
\omega_{k}=&\frac{1}{k\psi_0}\sum_{j=1}^{3}[\phi_{j-1}-(k-j)\psi_{j}]\omega_{k-j}, \quad k \geq 3,
\end{split}\end{equation}
where,
\begin{equation}\label{NC.2.2}\begin{split}
\phi_0=&2\alpha(\theta-1)(\alpha\theta+1)+\alpha\theta(\theta-\frac{3}{2}),
\\
\phi_1=&-\alpha(2\theta^2-3\alpha\theta+4\alpha\theta^2-1),
\\
\phi_2=&-\alpha\theta(\frac{1}{2}-\theta+\alpha-2\alpha\theta),
\end{split}\end{equation}
and
\begin{equation}\label{NC.2.3}\begin{split}
\psi_0=&\frac{1}{2}(3-2\theta)(1+\alpha\theta),
\quad\psi_1=-\frac{\alpha\theta}{2}(3-2\theta)-2(1-\theta)(\alpha\theta+1),
\\
\psi_2=&-\frac{1}{2}(\alpha\theta+1)(2\theta-1)-2\alpha\theta(\theta-1),
\quad\psi_3=-\frac{1}{2}\alpha\theta(1-2\theta).
\end{split}\end{equation}
\end{lemma}
With the analysis in \cite{LiuYin} we have the estimate that, if $\varphi(t)=t^{\beta}$ with $\beta>0$, then
\begin{equation}\label{N.3}
{}_{RL} D_{0,t}^{\alpha}\varphi(t_n)=D_{\tau,\omega}^{\alpha}\varphi^n+E^n,
\quad \text{for }\alpha \in (0,1],
\end{equation}
where $E^n=O(\tau^2)$.
Note that when $\alpha=1$, the operator ${}_{RL} D_{0,t}^{\alpha}$ is defined as the traditional first derivative.
We take the convolution weights $\omega_k^{(1)}$ with the assumption $\theta=0$, i.e., $u_t$ is approximated by traditional BDF2.
\par
Considering $u_0 \neq 0$, we take $v=u-u_0$ in which case equation (\ref{I.1}) can be formulated as
\begin{equation}\label{N.4}
v_t={}_{RL}D_{0,t}^{1-\gamma}\Delta v-\mu^2{}_{RL}D_{0,t}^{1-\kappa}v+F(\boldsymbol z,t),
\end{equation}
where $F(\boldsymbol z,t)=f(\boldsymbol z,t)+\Delta u_0 \frac{t^{\gamma-1}}{\Gamma(\gamma)}-\mu^2 u_0 \frac{t^{\kappa-1}}{\Gamma(\kappa)}$.
With the relation (\ref{N.3}), we can get
\begin{equation}\label{N.5}
D_{\tau,\omega}^{1}v^n=D_{\tau,\omega}^{1-\gamma}\Delta v^n-\mu^2D_{\tau,\omega}^{1-\kappa}v^n+F^n+E^n,
\end{equation}
where $F^n=F(\boldsymbol z,t_n)$.
\par
To derive the fully discrete scheme, we define $X_h$ as the subspace of $H_0^1(\Omega)$ as follows
\begin{equation}\label{N.6}\begin{split}
X_h=\{\chi \in H^{1}_0(\Omega):\chi|_{e} \in \mathbb{P}_r(\boldsymbol z), e \in \mathcal{T}_h\},
\end{split}
\end{equation}
where $\mathbb{P}_r(\boldsymbol z)$ is the set of linear polynomials of $\boldsymbol z$ with the degree no greater than $r$$(r \in \mathbb{Z}^{+})$ in one variable.
Denote $\mathcal{T}_h$ as a shape-regular and quasi-uniform triangulation of $\Omega$, and denote by $h$ the mesh size of $\mathcal{T}_h$.
Then the fully discrete scheme of equation (\ref{N.4}) is to find $V^n:[0,T] \longmapsto X_h$, such that
\begin{equation}\label{N.7}
D_{\tau,\omega}^{1}(V^n,\chi_h)+D_{\tau,\omega}^{1-\gamma}(\nabla V^n,\nabla \chi_h)+\mu^2 D_{\tau,\omega}^{1-\kappa}(V^n,\chi_h)=(F^n,\chi_h),
\end{equation}
holds for any $\chi_h \in X_h$.
\section{Stability analysis}
In this section we derive the stability estimate for the fully discrete scheme (\ref{N.7}). Considering the starting part dose not affect the stability, we next mainly analyse the following scheme, after omitting the starting part from (\ref{N.7}),
\begin{equation}\label{S.1}\begin{split}
\tau^{-1}\sum_{k=1}^{n}\omega^{(1)}_{n-k}(V^k,\chi_h)
+\tau^{\gamma-1}\sum_{k=1}^{n}\omega^{(1-\gamma)}_{n-k}(\nabla V^k,\nabla \chi_h)
\\
+\mu^2 \tau^{\kappa-1}\sum_{k=1}^{n}\omega^{(1-\kappa)}_{n-k}(V^k,\chi_h)=(F^n,\chi_h),
\end{split}\end{equation}
for any $\chi_h \in X_h$. The index $k$ starts from $1$ is due to the fact that $V^0=0$.
\par
First, we introduce some lemmas about the fractional $\theta$-methods which are crucial for the stability analysis.
\begin{lemma}(\textbf{Szego's theorem})(See \cite{Fisher,Szego})\label{l.0}
If the generating function $G(x)=\sum_{k=-\infty}^{\infty}c_k e^{ikx}$ of a symmetric Toeplitz matrix $D_n$ is a (almost everywhere existing) derivative of a real monotonically nondecreasing function, then
\begin{equation}\label{A.1}\begin{split}
\lim_{n\to\infty}\det(D_n)/\det(D_{n-1})=\exp\bigg(\frac{1}{2\pi}\int_{0}^{2\pi}\ln G(x)\mathrm{d}x\bigg),
\end{split}
\end{equation}
where $D_n$ is defined as
\[D_{n}=
\begin{pmatrix}
  c_0 & c_1 & c_2 & \cdots & c_n \\
  c_{-1} & c_0 & c_1 & \cdots & c_{n-1} \\
  c_{-2} & c_{-1} & c_0 & \cdots & c_{n-2} \\
  \vdots & \vdots & \vdots & \ddots & \vdots \\
  c_{-n} & c_{-(n-1)} & c_{-(n-2)} & \cdots & c_0
\end{pmatrix},\quad
\text{with }c_{-k}=c_k,
\]
and the limit in (\ref{A.1}) is approached from above.
\end{lemma}
\begin{lemma}\label{l.1}
Assume $\alpha \in (0,1]$, and the sequence $\{\omega_k^{(\alpha)}\}$ is generated by  (\ref{N.2}) for the FBT-$\theta$ or FBN-$\theta$ method. For any vector $(v^0,\cdots,v^{n-1})\in \mathbb{R}^{n}$ with $n\geq 1$, we have the following estimate
\begin{equation}\label{S.2}\begin{split}
\sum_{j=0}^{n-1}v^j\sum_{k=0}^{j}\omega^{(\alpha)}_{j-k} v^{k} \geq 0.
\end{split}
\end{equation}
Furthermore, when $\alpha=1$, the inequality (\ref{S.2}) can be strengthened as
\begin{equation}\label{S.3}\begin{split}
\sum_{j=0}^{n-1}v^j\sum_{k=0}^{j}\omega^{(1)}_{j-k} v^{k} \geq \varepsilon_0 (v^{n-1})^2,
\end{split}
\end{equation}
where the constant $\varepsilon_0$ is positive and independent of $n$ and the vector \\$(v^0,\cdots,v^{n-1})$.
\end{lemma}
\textbf{Proof. }
Let $c_0=\omega^{(\alpha)}_0$ and $c_{-k}=c_k=\omega_k^{(\alpha)}/2~(k=1,2,\cdots)$.
The left hand side of (\ref{S.2}) can be formulated as $\sum c_{j-k}v^j v^{k} (j,k=0,1,\cdots,n-1)$, which is the Toeplitz form (see \cite{Szego}) associated with the generating function
\begin{equation}\label{S.4}\begin{split}
f_{\alpha}(x)=&\sum_{k=-\infty}^{\infty}c_k e^{ikx}
=\omega_0^{(\alpha)}+\frac{1}{2}\sum_{k=1}^{\infty}\omega^{(\alpha)}_k e^{ikx}
+\frac{1}{2}\sum_{k=1}^{\infty}\omega^{(\alpha)}_k e^{-ikx}
\\=&
\frac{1}{2}\omega^{(\alpha)}(e^{ix})+\frac{1}{2}\omega^{(\alpha)}(e^{-ix}),
\quad x\in[0,2\pi].
\end{split}
\end{equation}
Considering the theorem on p.19 \cite{Szego} and the fact that $f_{\alpha}(x)$ is symmetric with respect to $x=\pi$, the inequality (\ref{S.2}) holds provided $f_{\alpha}(x)$ is nonnegative for $x\in [0,\pi]$ with fixed $\alpha \in (0,1]$. Actually, for the FBT-$\theta$ method, we have
\begin{equation}\label{S.5}\begin{split}
\omega^{(\alpha)}(\xi)=\bigg(\frac{3-2\theta}{2-2\theta}\bigg)^{\alpha}(1-\xi)^{\alpha}(1-\lambda_1\xi)^{-\alpha}(1-\lambda_2\xi)^{\alpha},
\end{split}
\end{equation}
where $\lambda_1=\frac{\theta}{\theta-1}$ and $\lambda_2=\frac{1-2\theta}{3-2\theta}$.
Note that $\lambda_1 \in (-1,1)$ and $\lambda_2 \in (0,1)$ since $\theta \in (-\infty,\frac{1}{2})$.
With the help of the equalities (see theorem 9 on p.78, \cite{LiChang})
\begin{equation}\label{S.6}\begin{split}
(1-e^{\pm ix})^{\alpha}=&\big(2\sin(x/2)\big)^{\alpha}e^{\pm\frac{ i\alpha}{2}(x-\pi)},
\\
(x-yi)^{\alpha}=&(x^2+y^2)^{\frac{\alpha}{2}}e^{i\alpha\phi}, ~\phi=-\arctan\frac{y}{x},
\end{split}
\end{equation}
combining (\ref{S.4}), we can get
\begin{equation}\label{S.7}\begin{split}
f_{\alpha}(x)=\bigg(\frac{3-2\theta}{2-2\theta}\bigg)^{\alpha}\bigg(2\sin\frac{x}{2}\bigg)^{\alpha}
\bigg(\frac{1+\lambda_2^2-2\lambda_2 \cos x}{1+\lambda_1^2-2\lambda_1\cos x}\bigg)^{\frac{\alpha}{2}}g_{\alpha,\theta}(x),
\end{split}
\end{equation}
where $g_{\alpha,\theta}(x)=\cos \alpha(\frac{x}{2}-\frac{\pi}{2}+\phi_2-\phi_1)$, and $\phi_i~(i=1,2)$ are defined as the following
\begin{equation}\label{S.8}\begin{split}
\phi_1=\arctan\frac{\lambda_1 \sin x}{\lambda_1\cos x-1},
\quad
\phi_2=\arctan\frac{\lambda_2 \sin x}{\lambda_2\cos x-1}.
\end{split}
\end{equation}
Let $h_{\theta}(x):=\frac{x}{2}-\frac{\pi}{2}+\phi_2-\phi_1$.
Next we show that $h_{\theta}(x) \in [-\frac{\pi}{2},0]$ for any $x\in[0,\pi]$ with fixed $\theta \in (-\infty,\frac{1}{2})$.
Take the first derivative of $h_{\theta}(x)$ to derive that
\begin{equation}\label{S.9}\begin{split}
h'_{\theta}(x)=\frac{\Lambda_{\theta}(\cos x)}{2(\lambda_1^2-2\lambda_1\cos x+1)(\lambda_2^2-2\lambda_2\cos x+1)},
\end{split}
\end{equation}
where $\Lambda_{\theta}(t)=4\lambda_1\lambda_2t^2-4\lambda_2(1+\lambda_1\lambda_2)t+3\lambda_2^2-\lambda_1^2+\lambda_1^2\lambda_2^2+1$.
Careful examination shows that the minimum of $\Lambda_{\theta}(t)$ for $t\in [-1,1]$ can be taken only at the end points of the interval $[-1,1]$.
It is easy to check that $\Lambda_{\theta}(1)=0$ and $\Lambda_{\theta}(1)>0$ and we omit the proof here.
Hence, $h_{\theta}(x)$ is a monotone nondecreasing function, and $h_{\theta}(x)\in [-\frac{\pi}{2},0]$, in which case $g_{\alpha,\theta}(x)\geq0$ and $f_{\alpha}(x)\geq 0$.
\par
\begin{figure}[h]
\begin{center}
\begin{minipage}{12cm}
  \centering\includegraphics[width=12cm]{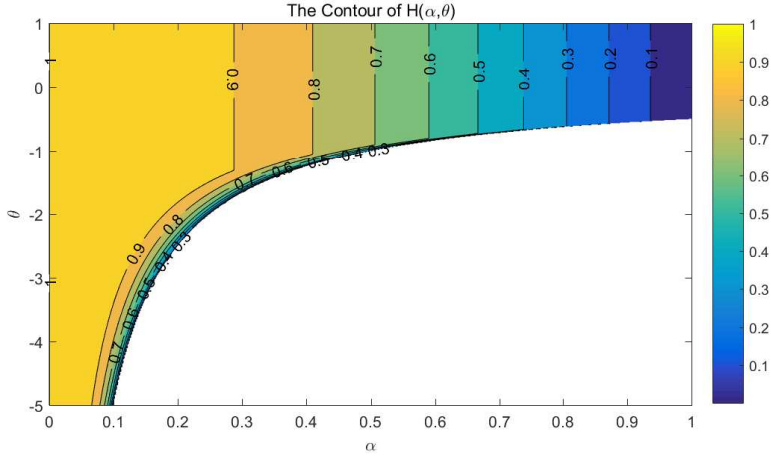}
  \caption{The contour of the function $H(\alpha,\theta)$ with $\alpha \in (0,1],\theta\in[-\frac{1}{2\alpha},1]$.}\label{C1}
\end{minipage}
\end{center}
\end{figure}

For the FBN-$\theta$ method, we have
\begin{equation}\label{S.10}\begin{split}
\omega^{(\alpha)}(\xi)=(\frac{3}{2}-\theta)^{\alpha}(1+\alpha\theta)(1-\xi)^{\alpha}(1-\lambda'_1\xi)(1-\lambda'_2\xi)^{\alpha},
\end{split}
\end{equation}
where $\lambda'_1=\frac{\alpha\theta}{1+\alpha\theta}$ and $\lambda'_2=\frac{1-2\theta}{3-2\theta}$.
Similar to the analysis of the FBT-$\theta$ method, we have
\begin{equation}\label{S.11}\begin{split}
f_{\alpha}(x)=(\frac{3}{2}-\theta)^{\alpha}(1+\alpha\theta)(2\sin\frac{x}{2})^{\alpha}
(1+\lambda^{'2}_1-2\lambda^{'}_1\cos x)^{\frac{1}{2}}
\\
(1+\lambda^{'2}_2-2\lambda^{'}_2\cos x)^{\frac{\alpha}{2}}g_{\alpha,\theta}(x),
\end{split}
\end{equation}
where $g_{\alpha,\theta}(x)=\cos\big(\frac{\alpha}{2}(x-\pi)+\phi_1+\alpha\phi_2\big)$ and $\phi_i~(i=1,2)$ are defined in (\ref{S.8}) with $\lambda_i$ replaced by $\lambda'_i$.
To analytically prove that $g_{\alpha,\theta}(x)$ is nonnegative is a tedious work, hence, here we merely numerically demonstrate $g_{\alpha,\theta}(x)\geq 0$ for $(x,\alpha,\theta)\in [0,\pi]\times\{(\alpha,\theta):0<\alpha\leq 1, -\frac{1}{2\alpha}\leq \theta \leq 1 \}$.
To this end, define the function $H(\alpha,\theta):=\min_{0\leq x\leq \pi}g_{\alpha,\theta}(x)$.
From Fig. \ref{C1}, one can easily check that $H(\alpha,\theta)$ is nonnegative, which means the function $g_{\alpha,\theta}(x)$ as well as $f_{\alpha}(x)$ is nonnegative.
\par
Moreover, by careful examination of the contours of $H(\alpha,\theta)$, one may find that for fixed $\alpha \in (0,1]$, the value of $H(\alpha,\theta)$ is not affected by $\theta$ so long as $\theta$ is far away from the curve depicted by the function $\theta=-\frac{1}{2\alpha}$.
\par
To prove the inequality (\ref{S.3}), we just need to prove that the matrix $D_n-\varepsilon_0E_n$ is positive semi-definite for any $n \geq 1$, where the entries $e_{ij}$ of $E_{n}$ are zeros except for $e_{nn}=1$, and $D_n$ is a Toeplitz matrix defined as
\[D_{n}=
\begin{pmatrix}
  c_0 & c_1 & c_2 & \cdots & c_n \\
  c_{-1} & c_0 & c_1 & \cdots & c_{n-1} \\
  c_{-2} & c_{-1} & c_0 & \cdots & c_{n-2} \\
  \vdots & \vdots & \vdots & \ddots & \vdots \\
  c_{-n} & c_{-(n-1)} & c_{-(n-2)} & \cdots & c_0
\end{pmatrix},\quad
\text{with }c_0=\omega^{(\alpha)}_0, ~c_{-k}=c_k=\frac{\omega_k^{(\alpha)}}{2}.
\]
Note that $\det(D_n-\varepsilon_0E_n)=\det(D_n)-\varepsilon_0\det(D_{n-1})$ and that a matrix is positive semi-definite if and only if all its principal minors are nonnegative, by Lemma \ref{l.0}, we take $\varepsilon_0$ satisfying
\begin{equation}\label{S.12}\begin{split}
\varepsilon_0 = \lim_{n\to\infty}\det(D_n)/\det(D_{n-1})=\exp\bigg(\frac{1}{2\pi}\int_{0}^{2\pi}\ln f_{\alpha}(x)\mathrm{d}x\bigg),
\end{split}
\end{equation}
where $f_{\alpha}(x)=f_1(x)=\sqrt{10-6\cos x}\sin \frac{x}{2} \sin(\frac{x}{2}+\arctan \frac{\sin x}{\cos x-3})=4(\sin \frac{x}{2})^4$ by letting $\theta=0$.
Considering the integral equality
\begin{equation}\label{S.13}\begin{split}
\int_{0}^{2\pi}\ln \sin \frac{x}{2}\mathrm{d}x=-2\pi \ln 2,
\end{split}
\end{equation}
we can easily derive that $\varepsilon_0=\frac{1}{4}$, which completes the proof of the lemma.
\begin{remark}
We remark that for the inequality (\ref{S.3}) similar result has been derived by Gao et al. \cite{Sun1}, however, the proof of (\ref{S.3}) in this paper is based on the theory of Toeplitz forms in which the generating function plays a crucial role.
The most important thing is that there always exists a positive constant $\varepsilon_0$ provided the generating function $f_{\alpha}(x)$ is the (almost everywhere existing) derivative of a real monotonically nondecreasing function and the integrand $\ln f_{\alpha}(x)$ is Lebesgue integrable (see theorem 1 on p.336, \cite{Fisher}).
\end{remark}
\begin{theorem}\label{T.1}
Suppose the solution $u(\cdot,t)$ is of the form (\ref{N.1.1}) and is sufficiently smooth, i.e., with $\sigma_1 \geq 3$. Let $V^n$ be the numerical solution of (\ref{N.4}), and define $U^n:=V^n+u_{0h}$ where $u_{0h} \in X_h$ is a proper approximation to $u_0$.
Then for sufficiently small $\tau$, the scheme (\ref{N.7}) is unconditionally stable with the following estimate
\begin{equation}\label{ST.1}\begin{split}
\|U^n\|\leq &
C\bigg(\tau\sum_{j=1}^{n}\|f^j\| +\|\Delta u_{0}\|+\|u_{0}\|+\|u_{0h}\|\bigg),
\end{split}
\end{equation}
where $C$ is independent of $n$, $h$ and $\tau$.
\end{theorem}
\textbf{Proof. }
We multiply both sides of (\ref{S.1}) by $\tau$, replace $n$ with $j$, take $\chi_h$ as $V^j$ and sum the index $j$ from $1$ to $n$ to formulate
\begin{equation}\label{ST.2}\begin{split}
\sum_{j=1}^{n}\sum_{k=1}^{j}\omega^{(1)}_{j-k}(V^k,V^j)
+\tau^{\gamma}\sum_{j=1}^{n}\sum_{k=1}^{j}\omega^{(1-\gamma)}_{j-k}(\nabla V^k,\nabla V^j)
\\
+\mu^2\tau^{\kappa}\sum_{j=1}^{n}\sum_{k=1}^{j}\omega^{(1-\kappa)}_{j-k}(V^k,V^j)
=\tau\sum_{j=1}^{n}(F^j,V^j).
\end{split}\end{equation}
Then by Lemma \ref{l.1}, Cauchy-Schwarz inequality and Young inequality, we have
\begin{equation}\label{ST.3}\begin{split}
\varepsilon_0\|V^n\|^2 \leq \varepsilon \max_{1\leq j \leq n}\|V^j\|^2+\frac{1}{4\varepsilon}\bigg(\tau\sum_{j=1}^{n}\|F^j\|\bigg)^2.
\end{split}\end{equation}
For the second term on right hand side of (\ref{ST.3}), we have the estimate that
\begin{equation}\label{ST.4}\begin{split}
\frac{1}{4\varepsilon}\bigg(\tau\sum_{j=1}^{n}\|F^j\|\bigg)^2
\leq &
\frac{1}{4\varepsilon}\bigg(\tau\sum_{j=1}^{n}\|f^j\|+\|\Delta u_{0}\| \tau\sum_{j=1}^{n}\frac{t_j^{\gamma-1}}{\Gamma(\gamma)}+\mu^2 \|u_{0}\|\tau\sum_{j=1}^{n} \frac{t_j^{\kappa-1}}{\Gamma(\kappa)}\bigg)^2
\\\leq &
\frac{1}{4\varepsilon}\bigg(\tau\sum_{j=1}^{n}\|f^j\|+ \frac{T^{\gamma}}{\Gamma(\gamma+1)}\|\Delta u_{0}\|+\frac{\mu^2T^{\gamma}}{\Gamma(\gamma+1)} \|u_{0}\|\bigg)^2
\\\leq &
C\bigg(\tau\sum_{j=1}^{n}\|f^j\|\bigg)^2
+C\|\Delta u_{0}\|^2+C\|u_{0}\|^2.
\end{split}\end{equation}
If we take $\varepsilon < \varepsilon_0$, considering (\ref{ST.3}) and (\ref{ST.4}), we can get
\begin{equation}\label{ST.5}\begin{split}
\|V^n\|^2\leq &
C\bigg(\tau\sum_{j=1}^{n}\|f^j\|\bigg)^2
+C\|\Delta u_{0}\|^2+C\|u_{0}\|^2.
\end{split}\end{equation}
Replacing $V^n$ by $U^n-u_{0h}$ in (\ref{ST.5}) and combining the triangle inequality, the proof of the theorem is completed.
\section{Error analysis}
In this section, we derive the optimal error estimate of the numerical scheme. Note that in our theoretical analysis, the starting part is omitted by assuming the solution is sufficiently smooth.
We first define the projection operator $R_h:H_0^1 \to X_h$ such that for given $u\in H_0^1(\Omega)$, it holds that
\begin{equation}\label{E.1}
(\nabla R_h u, \nabla \chi)=(\nabla u,\nabla \chi),\quad \forall \chi \in X_h.
\end{equation}
For the operator $R_h$, we assume the following estimate (see \cite{Vidar})
\begin{equation}\label{E.2}
\|R_h u-u\|+h\|\nabla R_h u-\nabla u\| \leq Ch^{r+1}\|u\|_{r+1}, \quad \text{for}~ u \in H^{r+1} \cap H_0^1,
\end{equation}
where $C$ is independent of $h$ and $u$.
\begin{theorem}\label{T.2}
Let $u$ be the solution of the equation (\ref{I.1}) of the form (\ref{N.1.1}), $v^n=u^n-u_0$ be the solution of (\ref{N.4}) and $V^n$ be the solution of (\ref{N.7}).
Define $U^n:=V^n+u_{0h}$ with $U^0=u_{0h}=R_h u_0$, then $U^n$ is the approximation of $u^n$. Suppose $u\in C([0,T]; H^{r+1}(\Omega) \cap H_0^1(\Omega)) \cap C^3([0,T]; L^2(\Omega))$, then we have the error estimate
\begin{equation}\label{4.3}
\|U^n-u(t_n)\| \leq C (\tau^2+h^{r+1})
\end{equation}
uniformly for $n=1,2,\cdots,N$, where the constant $C$ is independent of $h$ and $\tau$.
\end{theorem}
\textbf{Proof.} Let $V^n-v(t_n)=(V^n-R_h v^n)+(R_h v^n-v^n)=:\eta^n + \rho^n$.
By integrating both sides of (\ref{N.5}) with $\chi_h$ on $\Omega$ and combining with the fully discrete scheme (\ref{N.7}) as well as (\ref{E.1}), we have
\begin{equation}\label{TE.2}\begin{split}
&D_{\tau,\omega}^{1}(\eta^n,\chi_h)
+D_{\tau,\omega}^{1-\gamma}(\nabla \eta^n,\nabla \chi_h)+\mu^2 D_{\tau,\omega}^{1-\kappa}(\eta^n,\chi_h)
\\
=&-(E^n,\chi_h)-D_{\tau,\omega}^{1}(\rho^n,\chi_h)
-\mu^2 D_{\tau,\omega}^{1-\kappa}(\rho^n,\chi_h).
\end{split}\end{equation}
Quite similar to the analysis in theorem \ref{T.2}, we multiply both sides of (\ref{TE.2}) by $\tau$, replace $n$ with $j$, take $\chi_h$ as $\eta^j$ and sum the index $j$ from $1$ to $n$ to obtain
\begin{equation}\label{TE.3}\begin{split}
\sum_{j=1}^{n}\sum_{k=1}^{j}\omega_{j-k}^{(1)}(\eta^k,\eta^j)
+\tau^{\gamma}\sum_{j=1}^{n}\sum_{k=1}^{j}\omega_{j-k}^{(1-\gamma)}(\nabla \eta^k,\nabla \eta^j)
\\+\mu^2\tau^{\kappa} \sum_{j=1}^{n}\sum_{k=1}^{j}\omega_{j-k}^{(1-\kappa)}(\eta^k,\eta^j)
\\
=-\tau\sum_{j=1}^{n}(E^j,\eta^j)
-\tau\sum_{j=1}^{n}(D_{\tau,\omega}^{1}\rho^j,\eta^j)
-\mu^2 \tau\sum_{j=1}^{n}(D_{\tau,\omega}^{1-\kappa}\rho^j,\eta^j).
\end{split}\end{equation}
The estimates for the right hand side of (\ref{TE.3}) are stated in the following.
\begin{equation}\label{TE.5}\begin{split}
\tau\sum_{j=1}^{n}(E^j,\eta^j)\leq \frac{T}{4\varepsilon}\max_{1\leq j\leq n}\|E^j\|^2+T\varepsilon\max_{1\leq j \leq n}\|\eta^j\|^2,
\end{split}\end{equation}
\begin{equation}\label{TE.6}\begin{split}
&\tau\sum_{j=1}^{n}(D_{\tau,\omega}^{1}\rho^j,\eta^j)
+\mu^2 \tau\sum_{j=1}^{n}(D_{\tau,\omega}^{1-\kappa}\rho^j,\eta^j)
\\\leq &
\tau\sum_{j=1}^{n}\|D_{\tau,\omega}^{1}\rho^j\|\|\eta^j\|
+\mu^2 \tau\sum_{j=1}^{n}\|D_{\tau,\omega}^{1-\kappa}\rho^j\|\|\eta^j\|
\\\leq &
\max_{1\leq j \leq n}\|\eta^j\|\bigg(\tau\sum_{j=1}^{n}\|D_{\tau,\omega}^1\rho^j-\rho_t(t_j)\|+\tau\sum_{j=1}^{n}\|\rho_t(t_j)\|\bigg)
\\&+
\mu^2\max_{1\leq j \leq n}\|\eta^j\|\bigg(\tau\sum_{j=1}^{n}\|D_{\tau,\omega}^{1-\kappa}\rho^j-{}_{RL}D_{0,t}^{1-\kappa}\rho(t_j)\|+\tau\sum_{j=1}^{n}\|{}_{RL}D_{0,t}^{1-\kappa}\rho(t_j)\|\bigg)
\\\leq &
\varepsilon(1+\mu^2)\max_{1\leq j \leq n}\|\eta^j\|^2
+\frac{CT^2(1+\mu^2)}{2\varepsilon}\tau^4
\\&+
\frac{C(1+\mu^2)}{2\varepsilon}\bigg[\bigg(\int_{0}^{T}\|v_t\|_{r+1}\mathrm{d}t\bigg)^2+\bigg(\int_{0}^{T}\|{}_{RL}D_{0,t}^{1-\kappa}v\|_{r+1}\mathrm{d}t\bigg)^2\bigg]h^{2r+2}.
\end{split}\end{equation}
Considering Lemma \ref{l.1} and taking $\varepsilon < \varepsilon_0\min\{1,\frac{1}{2T},\frac{1}{2(1+\mu^2)}\}$, we combine the estimates (\ref{TE.3})-(\ref{TE.6}) to get
\begin{equation}\label{TE.7}\begin{split}
\|\eta^n\|\leq C(\tau^2+h^{r+1}).
\end{split}\end{equation}
Finally, by (\ref{E.2}) and $\|U^n-u(t_n)\|\leq \|V^n-v(t_n)\|+\|R_h u_0-u_0\|$, we complete the proof of the theorem.
\section{Numerical tests}
In this section we conduct some numerical experiments to further confirm our theoretical analysis.
The error $E(\tau,h)=\max_{0 \leq n \leq N}\|u^n-U^n\|$ is recorded and the convergence rate are derived by the formulas
\begin{equation}\label{Nt.1}\begin{split}
\text{temporal order }=\log_2\frac{E(2\tau,h)}{E(\tau,h)},
\quad
\text{spatial order }=\log_2\frac{E(\tau,2h)}{E(\tau,h)}.
\end{split}\end{equation}
To overcome the singularity at initial value for solutions with weak regularity, we take the approximation formula (\ref{N.1}) with the starting part, and compare the results with those obtained without the starting part.
Hence, for clarity, we denote by $E_c(\tau,h)$ the error derived by the fractional $\theta$-methods with the starting part; and by $E_o(\tau,h)$ the error derived without the starting part.
Note that there exist two fractional derivatives in equation (\ref{I.1}), and we choose different parameter $\theta$ which are subscripted as $\theta_{\gamma}$ and $\theta_{\kappa}$ for the fractional $\theta$-methods to approximate the two fractional derivatives ${}_{RL}D_{0,t}^{1-\gamma}$ and ${}_{RL}D_{0,t}^{1-\kappa}$, respectively.
\subsection{Example of one-dimensional space}
We take $\Omega=(0,1)$, $\mu=1$ and $T=1$.
The interval $\Omega$ is divided into a uniform partition as $0=x_0<x_1<\cdots<x_{N_s}=1$ with $N_s>0$.
Let $h=1/N_s$.
Define by $X_h$ the space of piecewise linear polynomials.
The exact solution is taken as $u(x,t)=(1+{\color{red}t^{\gamma}}+{\color{red}t^{\kappa}}+t^3)\sin(2\pi x)$ which is of weak singularity at initial value.
The term $f$ can be derived by substituting the $u(x,t)$ into the equation (\ref{I.1}) and the expression is omitted here.
\par
In Table \ref{C1}, we choose different pairs of $(\gamma,\kappa)$ and for each pair we take different FBT-$\theta$ formulas by varying $\theta$ to approximate the fractional derivatives of the equation (\ref{I.1}) under fixed fine space mesh $h=\frac{1}{5000}$.
With the time mesh taken as $\tau=\frac{1}{10}, \frac{1}{20}, \frac{1}{40}, \frac{1}{80}$, respectively, one can see that the error $E_c(\tau,h)$ is generally smaller than $E_o(\tau,h)$ and the rate of $E_c(\tau,h)$ is of $2$ compared with that of $E_o(\tau,h)$, which is much smaller than the optimal convergence rate in time.
\par
\begin{figure}[h]
\begin{center}
\begin{minipage}{5.5cm}
  \centering\includegraphics[width=6cm]{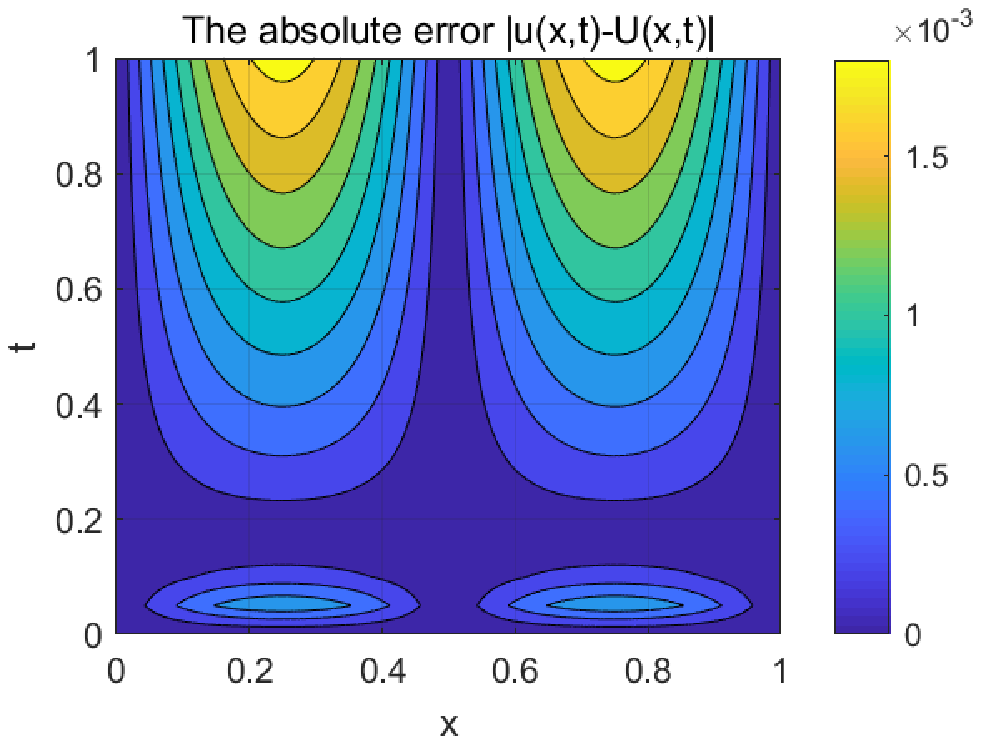}
  \caption{$h=\frac{1}{5000}$, $\tau=\frac{1}{20}$, $\gamma=0.6$, $\kappa=0.5$, $\theta_{\gamma}=0$, $\theta_{\kappa}=0.49$.}\label{C2}
\end{minipage}
\begin{minipage}{5.5cm}
  \centering\includegraphics[width=6cm]{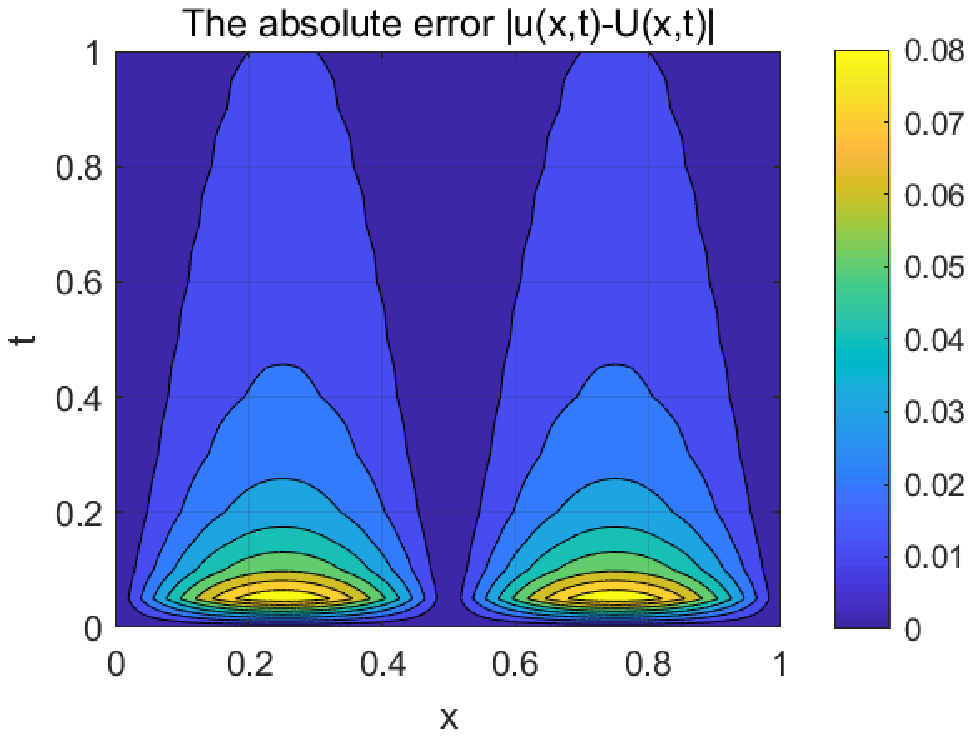}
  \caption{$h=\frac{1}{5000}$, $\tau=\frac{1}{20}$, $\gamma=0.6$, $\kappa=0.5$, $\theta_{\gamma}=0$, $\theta_{\kappa}=0.49$.}\label{C3}
\end{minipage}
\end{center}
\end{figure}

In Table \ref{C2}, we collect the errors and convergence rates in time of the experiment when using the FBN-$\theta$ method.
Note that $(\theta_{\gamma},\theta_{\kappa})$ for this method satisfies $\theta_{\gamma}\in [-\frac{1}{2(1-\gamma)},1]$ and $\theta_{\kappa}\in [-\frac{1}{2(1-\kappa)},1]$.
The fine space mesh is set as $h=\frac{1}{5000}$ and the time mesh is taken as $\tau=\frac{1}{10}, \frac{1}{20}, \frac{1}{40}, \frac{1}{80}$, respectively.
One can easily find out that the convergence rate in time is $2$ provided the starting part is added.
\begin{table}[h]
\centering
 \caption{The temporal convergence rate for the FBT-$\theta$ method with $h=\frac{1}{5000}$}\label{tab1}
 \begin{tabular}{ccccccc}
\toprule
  $(\gamma,\kappa)$ &$(\theta_{\gamma},\theta_{\kappa})$&  $\tau$ & $E_c(\tau, h)$ & rate &  $E_o(\tau, h)$ & rate \\
\midrule
	&	(0,0)	&	  1/10  	&	1.05368E-02	&	--	&	2.12608E-01	&	--	\\
	&		&	  1/20  	&	1.89214E-03	&	2.4773 	&	1.68531E-01	&	0.3352 	\\
	&		&	  1/40  	&	5.54574E-04	&	1.7706 	&	1.34676E-01	&	0.3235 	\\
	&		&	  1/80  	&	1.47119E-04	&	1.9144 	&	1.08343E-01	&	0.3139 	\\
(0.3,0.9)&	(0,0.49)	&	  1/10  	&	1.05368E-02	&	--	&	2.12639E-01	&	--	\\
	&		&	  1/20  	&	1.88811E-03	&	2.4804 	&	1.68546E-01	&	0.3353 	\\
	&		&	  1/40  	&	5.53517E-04	&	1.7702 	&	1.34683E-01	&	0.3236 	\\
	&		&	  1/80  	&	1.46853E-04	&	1.9143 	&	1.08346E-01	&	0.3139 	\\
	&	(-0.5,0.4)	&	  1/10  	&	1.05368E-02	&	--	&	1.93626E-01	&	--	\\
	&		&	  1/20  	&	3.48719E-03	&	1.5953 	&	1.54923E-01	&	0.3217 	\\
	&		&	  1/40  	&	9.75418E-04	&	1.8380 	&	1.24618E-01	&	0.3140 	\\
	&		&	  1/80  	&	2.55089E-04	&	1.9350 	&	1.00755E-01	&	0.3067 	\\
\midrule
	&	(-1,0.49)	&	  1/10  	&	1.07168E-02	&	--	&	6.60333E-02	&	--	\\
	&		&	  1/20  	&	3.08842E-03	&	1.7949 	&	4.90893E-02	&	0.4278 	\\
	&		&	  1/40  	&	8.19626E-04	&	1.9138 	&	3.70115E-02	&	0.4074 	\\
	&		&	  1/80  	&	2.10349E-04	&	1.9622 	&	2.83526E-02	&	0.3845 	\\
(0.6,0.5)&	(0.4,-1)	&	  1/10  	&	3.79088E-03	&	--	&	1.06932E-01	&	--	\\
	&		&	  1/20  	&	6.62162E-04	&	2.5173 	&	7.55151E-02	&	0.5019 	\\
	&		&	  1/40  	&	1.77051E-04	&	1.9030 	&	5.37684E-02	&	0.4900 	\\
	&		&	  1/80  	&	4.51734E-05	&	1.9706 	&	3.86893E-02	&	0.4748 	\\
	&	(-0.5,0)	&	  1/10  	&	8.02601E-03	&	--	&	7.20953E-02	&	--	\\
	&		&	  1/20  	&	2.26675E-03	&	1.8241 	&	5.29315E-02	&	0.4458 	\\
	&		&	  1/40  	&	5.95408E-04	&	1.9287 	&	3.93943E-02	&	0.4261 	\\
	&		&	  1/80  	&	1.51931E-04	&	1.9705 	&	2.97853E-02	&	0.4034 	\\
\midrule
	&	(0.49,0.49)	&	  1/10  	&	2.63144E-03	&	--	&	2.12803E-01	&	--	\\
	&		&	  1/20  	&	6.56725E-04	&	2.0025 	&	2.52160E-01	&	-0.2448 	\\
	&		&	  1/40  	&	1.75731E-04	&	1.9019 	&	2.88481E-01	&	-0.1941 	\\
	&		&	  1/80  	&	4.60032E-05	&	1.9336 	&	3.13814E-01	&	-0.1214 	\\
(0.9,0.1)&	(-0.1,0.49)	&	  1/10  	&	2.63144E-03	&	--	&	2.03506E-01	&	--	\\
	&		&	  1/20  	&	6.56725E-04	&	2.0025 	&	2.46765E-01	&	-0.2781 	\\
	&		&	  1/40  	&	1.75731E-04	&	1.9019 	&	2.85689E-01	&	-0.2113 	\\
	&		&	  1/80  	&	4.60032E-05	&	1.9336 	&	3.12513E-01	&	-0.1295 	\\
	&	(0.49,-1)	&	  1/10  	&	2.63144E-03	&	--	&	1.82841E-01	&	--	\\
	&		&	  1/20  	&	6.56725E-04	&	2.0025 	&	2.20926E-01	&	-0.2730 	\\
	&		&	  1/40  	&	1.75731E-04	&	1.9019 	&	2.59419E-01	&	-0.2317 	\\
	&		&	  1/80  	&	4.60032E-05	&	1.9336 	&	2.89589E-01	&	-0.1587 	\\
 \bottomrule
 \end{tabular}
\end{table}
\begin{table}[h]
\centering
 \caption{The temporal convergence rate for the FBN-$\theta$ method with $h=\frac{1}{5000}$}\label{tab2}
 \begin{tabular}{ccccccc}
\toprule
  $(\gamma,\kappa)$ &$(\theta_{\gamma},\theta_{\kappa})$&  $\tau$ & $E_c(\tau, h)$ & rate &  $E_o(\tau, h)$ & rate \\
\midrule
	&	(0,0)	&	  1/10  	&	6.80886E-03	&	--	&	1.37556E-01	&	--	\\
	&		&	  1/20  	&	1.83124E-03	&	1.8946 	&	1.01102E-01	&	0.4442 	\\
	&		&	  1/40  	&	5.00615E-04	&	1.8711 	&	7.49161E-02	&	0.4325 	\\
	&		&	  1/80  	&	1.29413E-04	&	1.9517 	&	5.59986E-02	&	0.4199 	\\
(0.4,0.8)	&	(0,0.5)	&	  1/10  	&	6.80886E-03	&	--	&	1.37586E-01	&	--	\\
	&		&	  1/20  	&	1.82753E-03	&	1.8975 	&	1.01118E-01	&	0.4443 	\\
	&		&	  1/40  	&	4.99636E-04	&	1.8709 	&	7.49240E-02	&	0.4325 	\\
	&		&	  1/80  	&	1.29161E-04	&	1.9517 	&	5.60028E-02	&	0.4199 	\\
	&	(0,1)	&	  1/10  	&	6.80886E-03	&	--	&	1.37480E-01	&	--	\\
	&		&	  1/20  	&	1.83526E-03	&	1.8914 	&	1.01063E-01	&	0.4440 	\\
	&		&	  1/40  	&	5.01635E-04	&	1.8713 	&	7.48956E-02	&	0.4323 	\\
	&		&	  1/80  	&	1.29670E-04	&	1.9518 	&	5.59879E-02	&	0.4198 	\\
\midrule
	&	(-1,-0.5)	&	  1/10  	&	2.60892E-02	&	--	&	1.50666E-01	&	--	\\
	&		&	  1/20  	&	7.38107E-03	&	1.8216 	&	9.67924E-02	&	0.6384 	\\
	&		&	  1/40  	&	1.94354E-03	&	1.9251 	&	6.02306E-02	&	0.6844 	\\
	&		&	  1/80  	&	4.97443E-04	&	1.9661 	&	3.68350E-02	&	0.7094 	\\
(0.5,0.6)	&	(-1,0.5)	&	  1/10  	&	2.59362E-02	&	--	&	1.49124E-01	&	--	\\
	&		&	  1/20  	&	7.33675E-03	&	1.8218 	&	9.59574E-02	&	0.6360 	\\
	&		&	  1/40  	&	1.93191E-03	&	1.9251 	&	5.98127E-02	&	0.6819 	\\
	&		&	  1/80  	&	4.94466E-04	&	1.9661 	&	3.66474E-02	&	0.7067 	\\
	&	(-1,1)	&	  1/10  	&	2.60129E-02	&	--	&	1.51346E-01	&	--	\\
	&		&	  1/20  	&	7.36178E-03	&	1.8211 	&	9.72006E-02	&	0.6388 	\\
	&		&	  1/40  	&	1.93842E-03	&	1.9252 	&	6.04715E-02	&	0.6847 	\\
	&		&	  1/80  	&	4.96117E-04	&	1.9661 	&	3.69737E-02	&	0.7098 	\\
\midrule
	&	(0.5,-0.5)	&	  1/10  	&	3.35206E-03	&	--	&	1.15488E-01	&	--	\\
	&		&	  1/20  	&	9.27758E-04	&	1.8532 	&	9.86337E-02	&	0.2276 	\\
	&		&	  1/40  	&	2.40945E-04	&	1.9450 	&	8.65064E-02	&	0.1893 	\\
	&		&	  1/80  	&	6.09667E-05	&	1.9826 	&	7.74200E-02	&	0.1601 	\\
(0.7,0.3)	&	(0.5,0.5)	&	  1/10  	&	3.35206E-03	&	--	&	1.20715E-01	&	--	\\
	&		&	  1/20  	&	8.32158E-04	&	2.0101 	&	1.03273E-01	&	0.2251 	\\
	&		&	  1/40  	&	2.16430E-04	&	1.9430 	&	9.05080E-02	&	0.1903 	\\
	&		&	  1/80  	&	5.47603E-05	&	1.9827 	&	8.07243E-02	&	0.1650 	\\
	&	(0.5,1)	&	  1/10  	&	3.35206E-03	&	--	&	1.15318E-01	&	--	\\
	&		&	  1/20  	&	9.07059E-04	&	1.8858 	&	9.84819E-02	&	0.2277 	\\
	&		&	  1/40  	&	2.35530E-04	&	1.9453 	&	8.63751E-02	&	0.1892 	\\
	&		&	  1/80  	&	5.95824E-05	&	1.9830 	&	7.73112E-02	&	0.1599 	\\
 \bottomrule
 \end{tabular}
\end{table}
\par
\begin{figure}[h]
\begin{center}
\begin{minipage}{5.5cm}
  \centering\includegraphics[width=6cm]{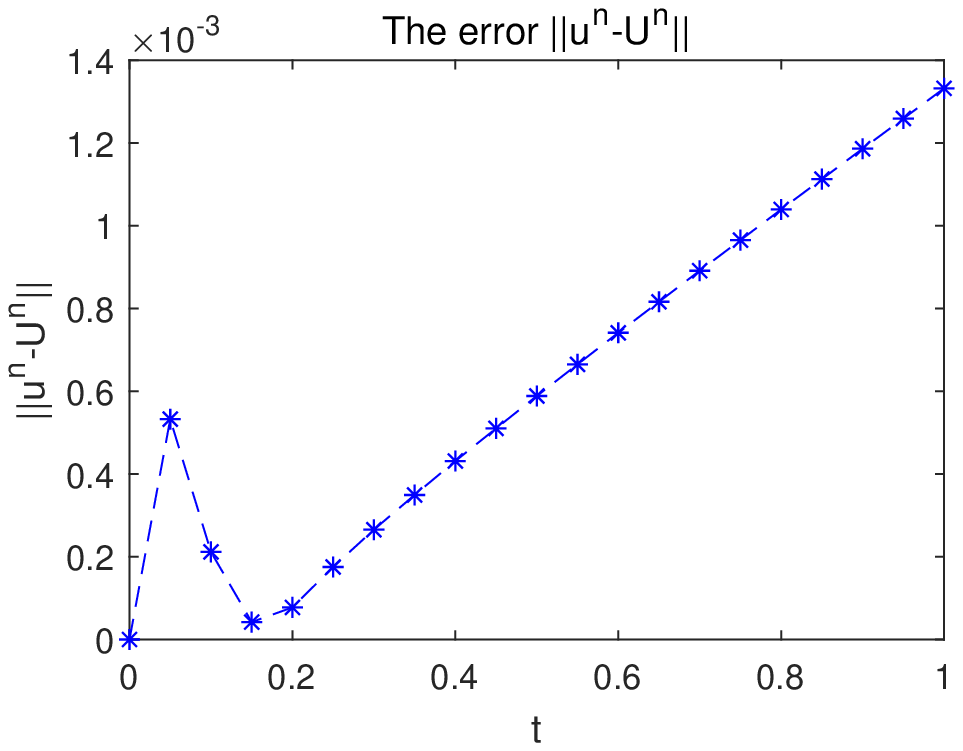}
  \caption{$h=\frac{1}{5000}$, $\tau=\frac{1}{20}$, $\gamma=0.6$, $\kappa=0.5$, $\theta_{\gamma}=0$, $\theta_{\kappa}=0.49$.}\label{C4}
\end{minipage}
\begin{minipage}{5.5cm}
  \centering\includegraphics[width=6cm]{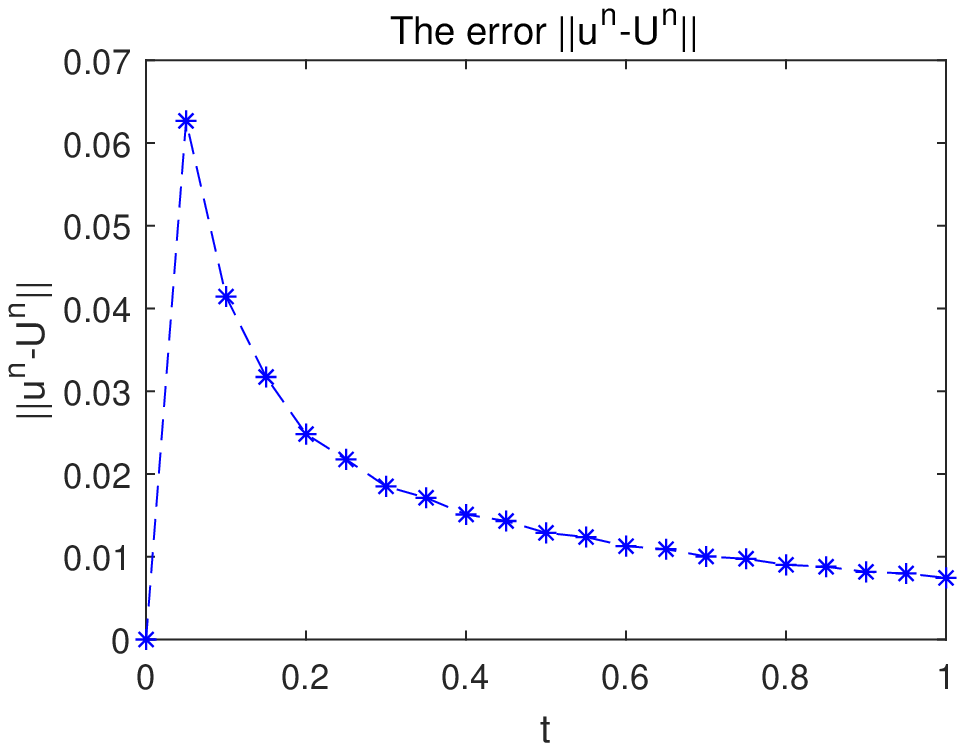}
  \caption{$h=\frac{1}{5000}$, $\tau=\frac{1}{20}$, $\gamma=0.6$, $\kappa=0.5$, $\theta_{\gamma}=0$, $\theta_{\kappa}=0.49$.}\label{C5}
\end{minipage}
\end{center}
\end{figure}

By a further examination of the error $|u-U|$ on the space-time plane when using the FBT-$\theta$ method, we find that with the starting part in our approximation formula, the error concentrates at the last time level (Fig. \ref{C2}), in contrast to the case with the starting part omitted, where the error mainly focuses on the several initial time levels (Fig. \ref{C3}).
To eliminate the effect of the space direction, we depict the error $\|u^n-U^n\|$ at each time level in Fig. \ref{C4} and Fig. \ref{C5} for the approximation formulas with and without the starting part, respectively.
A direct conclusion is that $E_c(\tau,h)$ is taken at the final time level, i.e., $\|u^N-U^N\|$, in contrast to $E_o(\tau,h)$ which is taken near initial time level.
\par
In Table \ref{tab3} and Table \ref{tab4}, we calculate the convergence order in space for both of the fractional $\theta$-methods with different pairs $(\theta_{\gamma},\theta_{\kappa})$.
With the fixed fine time mesh $\tau=\frac{1}{1000}$, we choose space meshes as $h=\frac{1}{10},\frac{1}{20},\frac{1}{40},\frac{1}{80}$, respectively, and record the errors $E_c(\tau,h)$.
One can see that the optimal second-order convergence rate is obtained, which confirms our theoretical results.
\begin{table}[h]
\centering
 \caption{The spatial convergence rate for the FBT-$\theta$ method with $\tau=\frac{1}{1000}$}\label{tab3}
 \begin{tabular}{ccccc}
\toprule
  $(\gamma,\kappa)$ &$(\theta_{\gamma},\theta_{\kappa})$&  $h$ & $E_c(\tau, h)$ & rate \\
\midrule
	&	(0,0)	&	  1/10  	&	9.73080E-02	&	--	\\
	&		&	  1/20  	&	2.44415E-02	&	1.9932 	\\
	&		&	  1/40  	&	6.11715E-03	&	1.9984 	\\
	&		&	  1/80  	&	1.52933E-03	&	2.0000 	\\
(0.6,0.2)	&	(0,0.4)	&	  1/10  	&	9.73080E-02	&	--	\\
	&		&	  1/20  	&	2.44415E-02	&	1.9932 	\\
	&		&	  1/40  	&	6.11717E-03	&	1.9984 	\\
	&		&	  1/80  	&	1.52934E-03	&	2.0000 	\\
	&	(-1,0.2)	&	  1/10  	&	9.73073E-02	&	--	\\
	&		&	  1/20  	&	2.44408E-02	&	1.9933 	\\
	&		&	  1/40  	&	6.11644E-03	&	1.9985 	\\
	&		&	  1/80  	&	1.52861E-03	&	2.0005 	\\
\bottomrule
 \end{tabular}
\end{table}
\begin{table}[h]
\centering
 \caption{The spatial convergence rate for the FBN-$\theta$ method with $\tau=\frac{1}{1000}$}\label{tab4}
 \begin{tabular}{ccccc}
\toprule
  $(\gamma,\kappa)$ &$(\theta_{\gamma},\theta_{\kappa})$&  $h$ & $E_c(\tau, h)$ & rate \\
\midrule
	&	(0,0)	&	  1/10  	&	9.67827E-02	&	--	\\
	&		&	  1/20  	&	2.43091E-02	&	1.9933 	\\
	&		&	  1/40  	&	6.08376E-03	&	1.9985 	\\
	&		&	  1/80  	&	1.52072E-03	&	2.0002 	\\
(0.3,0.9)	&	(1,0.5)	&	  1/10  	&	9.67819E-02	&	--	\\
	&		&	  1/20  	&	2.43082E-02	&	1.9933 	\\
	&		&	  1/40  	&	6.08286E-03	&	1.9986 	\\
	&		&	  1/80  	&	1.51982E-03	&	2.0008 	\\
	&	(-0.5,-1)	&	  1/10  	&	9.67816E-02	&	--	\\
	&		&	  1/20  	&	2.43079E-02	&	1.9933 	\\
	&		&	  1/40  	&	6.08257E-03	&	1.9987 	\\
	&		&	  1/80  	&	1.51953E-03	&	2.0011 	\\
\bottomrule
 \end{tabular}
\end{table}
\par
To further prove the necessity of adding the starting part, we next consider an example with the zero source term whose solution can not be expressed in a closed form.
In order to avoid too much complexity in calculating the exact solution, we take $\mu=0$ and $\Omega=(0,\pi)$ in which case the solution $u$ is
\begin{equation}\label{Nt.2}
u(x,t)=\sum_{j=0}^{\infty}\frac{(-t^\gamma)^j}{\Gamma(\gamma j+1)}\sin x,
\quad
u_0=\sin x.
\end{equation}
The temporal convergence rates are reported in Table \ref{tab4.1} and Table \ref{tab4.2} for the FBT-$\theta$ and FBN-$\theta$ methods, respectively, where $\gamma=0.8$ and $h=\frac{1}{5000}$.
From the column $E_c(\tau, h)$, one observes that the starting part can improve the scheme accuracy with which the optimal convergence rate is arrived at.
However, without the starting part the error is larger and the convergence rate is much lower as reported in the column $E_o(\tau, h)$ and the next column.

\begin{table}[h]
\centering
 \caption{The temporal convergence rate for the FBT-$\theta$ method with $h=\frac{1}{5000}$ and $\gamma=0.8$}\label{tab4.1}
 \begin{tabular}{cccccc}
\toprule
  $\theta$&  $\tau$ & $E_c(\tau, h)$ & rate & $E_o(\tau, h)$ & rate\\
\midrule
0	&	  1/20  	&	1.86654E-04	&	--	&	2.64217E-02	&	--	\\
	&	  1/40  	&	5.47490E-05	&	1.7695 	&	1.51117E-02	&	0.8061 	\\
	&	  1/80  	&	1.50683E-05	&	1.8613 	&	8.62130E-03	&	0.8097 	\\
	&	  1/160 	&	3.98128E-06	&	1.9202 	&	4.92431E-03	&	0.8080 	\\
0.49	&	  1/20  	&	2.24628E-04	&	--	&	2.68005E-02	&	--	\\
	&	  1/40  	&	6.55487E-05	&	1.7769 	&	1.52410E-02	&	0.8143 	\\
	&	  1/80  	&	1.79678E-05	&	1.8672 	&	8.66757E-03	&	0.8143 	\\
	&	  1/160 	&	4.73584E-06	&	1.9237 	&	4.94031E-03	&	0.8110 	\\
-0.5	&	  1/20  	&	1.48292E-04	&	--	&	2.61591E-02	&	--	\\
	&	  1/40  	&	4.36901E-05	&	1.7631 	&	1.50155E-02	&	0.8009 	\\
	&	  1/80  	&	1.20943E-05	&	1.8530 	&	8.58756E-03	&	0.8061 	\\
	&	  1/160 	&	3.20788E-06	&	1.9146 	&	4.91277E-03	&	0.8057 	\\
\bottomrule
 \end{tabular}
\end{table}
\begin{table}[h]
\centering
 \caption{The temporal convergence rate for the FBN-$\theta$ method with $h=\frac{1}{5000}$ and $\gamma=0.8$}\label{tab4.2}
 \begin{tabular}{cccccc}
\toprule
  $\theta$&  $\tau$ & $E_c(\tau, h)$ & rate & $E_o(\tau, h)$ & rate\\
\midrule
0	&	  1/20  	&	1.86654E-04	&	--	&	2.64217E-02	&	--	\\
	&	  1/40  	&	5.47490E-05	&	1.7695 	&	1.51117E-02	&	0.8061 	\\
	&	  1/80  	&	1.50683E-05	&	1.8613 	&	8.62130E-03	&	0.8097 	\\
	&	  1/160 	&	3.98128E-06	&	1.9202 	&	4.92431E-03	&	0.8080 	\\
0.5	&	  1/20  	&	2.02209E-04	&	--	&	2.65577E-02	&	--	\\
	&	  1/40  	&	5.91842E-05	&	1.7726 	&	1.51619E-02	&	0.8087 	\\
	&	  1/80  	&	1.62578E-05	&	1.8641 	&	8.63900E-03	&	0.8115 	\\
	&	  1/160 	&	4.29027E-06	&	1.9220 	&	4.93038E-03	&	0.8092 	\\
1	&	  1/20  	&	1.71142E-04	&	--	&	2.62564E-02	&	--	\\
	&	  1/40  	&	5.03715E-05	&	1.7645 	&	1.50475E-02	&	0.8031 	\\
	&	  1/80  	&	1.38928E-05	&	1.8583 	&	8.59807E-03	&	0.8074 	\\
	&	  1/160 	&	3.67493E-06	&	1.9185 	&	4.91622E-03	&	0.8065 	\\

\bottomrule
 \end{tabular}
\end{table}
\subsection{Example of two-dimensional space}
\begin{figure}[h]
\begin{center}
\begin{minipage}{5.5cm}
  \centering\includegraphics[width=6cm]{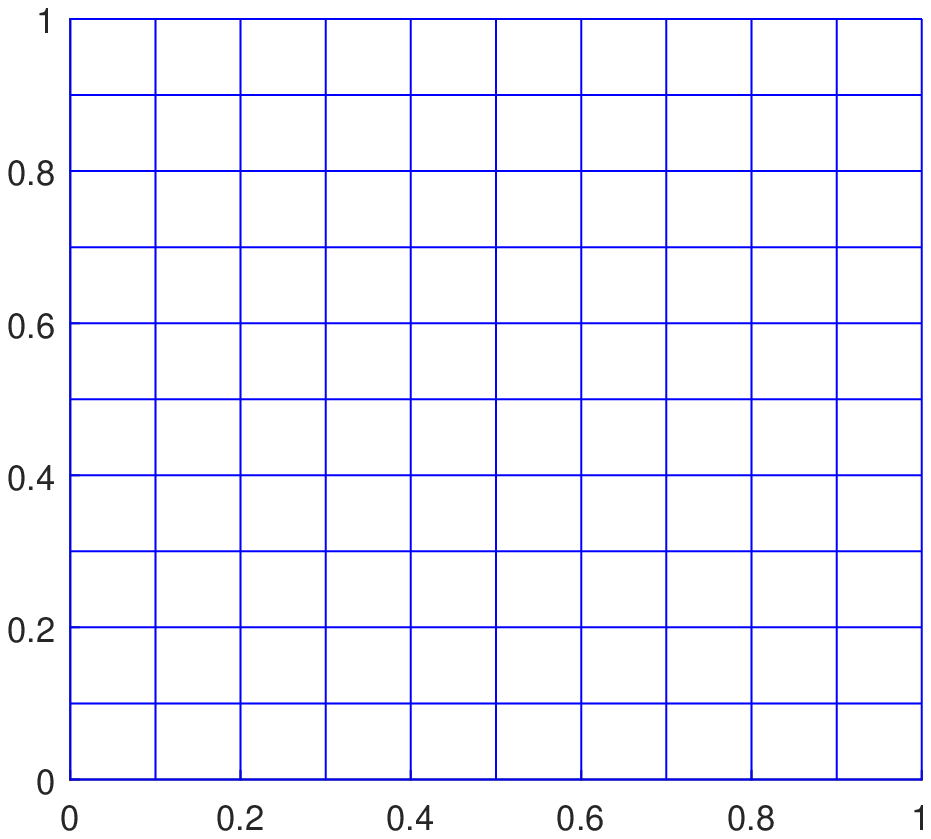}
  \caption{The space grid mesh with $h=\frac{\sqrt{2}}{10}$.}\label{C100}
\end{minipage}
\begin{minipage}{5.5cm}
  \centering\includegraphics[width=6cm]{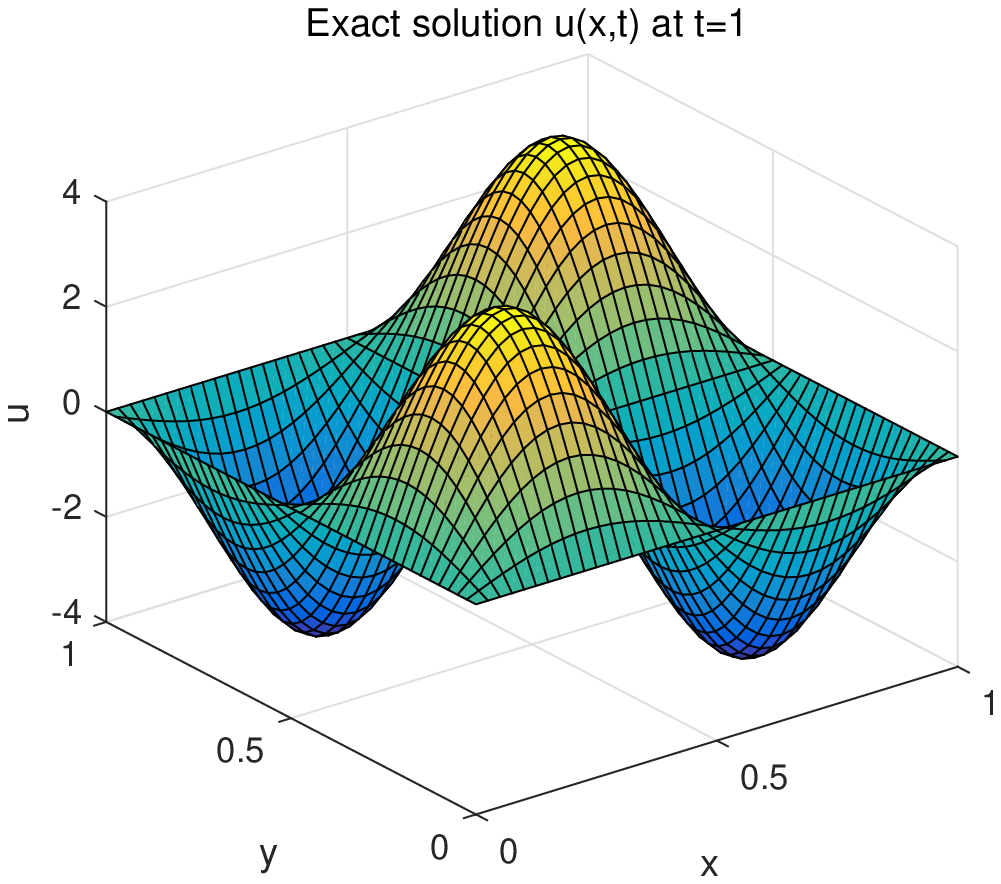}
  \caption{$h=\frac{\sqrt{2}}{40}$, $\tau=\frac{1}{200}$.}\label{C8}
\end{minipage}
\end{center}
\end{figure}

For the example of two-dimensional space, we take $\Omega=(0,1)\times(0,1)$, $T=1$ and $\mu=1$.
The triangulation $\mathcal{T}_h$ of $\Omega$ is as in Fig. \ref{C100}.
The elements of finite element space $X_h$ are chosen as the piecewise bilinear elements with the shape function $u=a_0xy+a_1x+a_2y+a_3$.
Take the exact solution of (\ref{I.1}) as $u(x,y,t)=(1+3t^3)\sin(2\pi x)\sin(2\pi y)$ and the source term $f$ can be obtained directly which is
\begin{equation*}\begin{split}
f(x,y,t)=&~\bigg[9t^2 + \frac{8t^{\gamma - 1}\pi^2(\gamma^3 + 3\gamma^2 + 2\gamma + 18t^3)}{\Gamma(\gamma + 3)}
\\&+ \frac{\mu^2 t^{\kappa - 1}(\kappa^3 + 3\kappa^2 + 2\kappa + 18t^3)}{\Gamma(\kappa + 3)}\bigg]\sin(2 \pi x)\sin(2 \pi y).
\end{split}\end{equation*}
\par
\begin{figure}[h]
\begin{center}
\begin{minipage}{5.5cm}
  \centering\includegraphics[width=6cm]{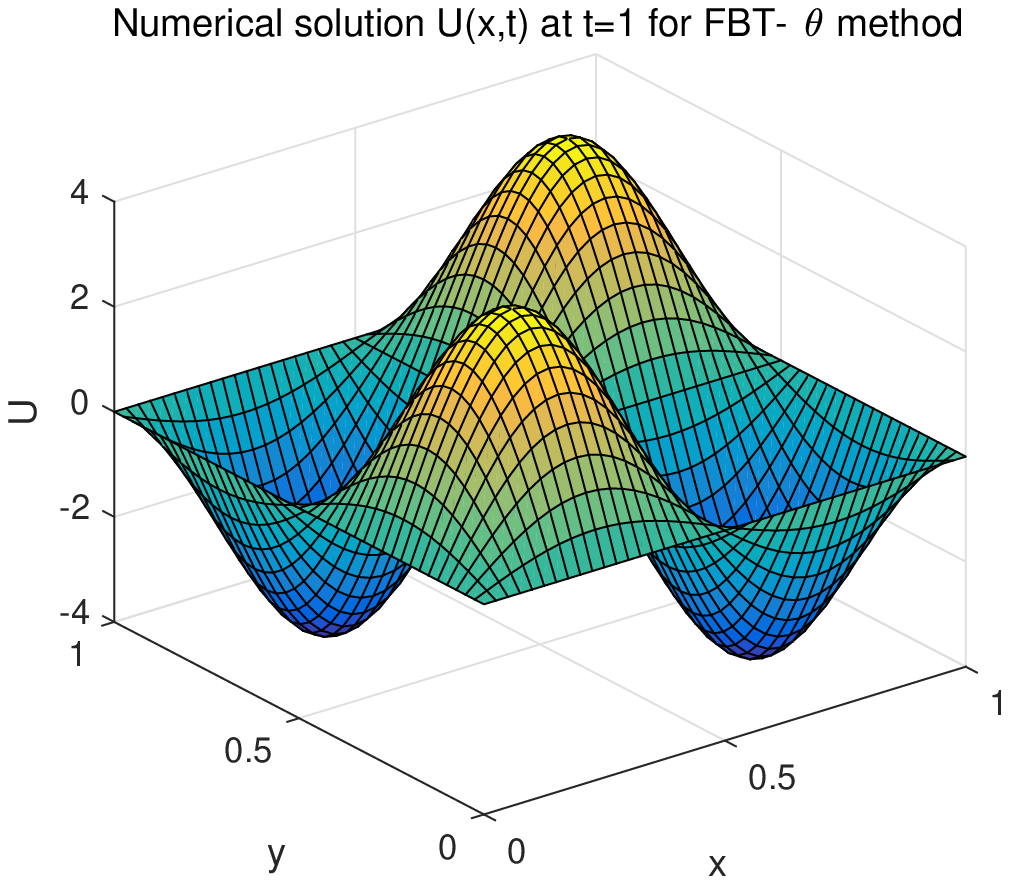}
  \caption{$h=\frac{\sqrt{2}}{40}$, $\tau=\frac{1}{200}$, $\gamma=0.4$, $\kappa=0.6$, $\theta_{\gamma}=0$, $\theta_{\kappa}=0.45$.}\label{C7}
\end{minipage}
\begin{minipage}{5.5cm}
  \centering\includegraphics[width=6cm]{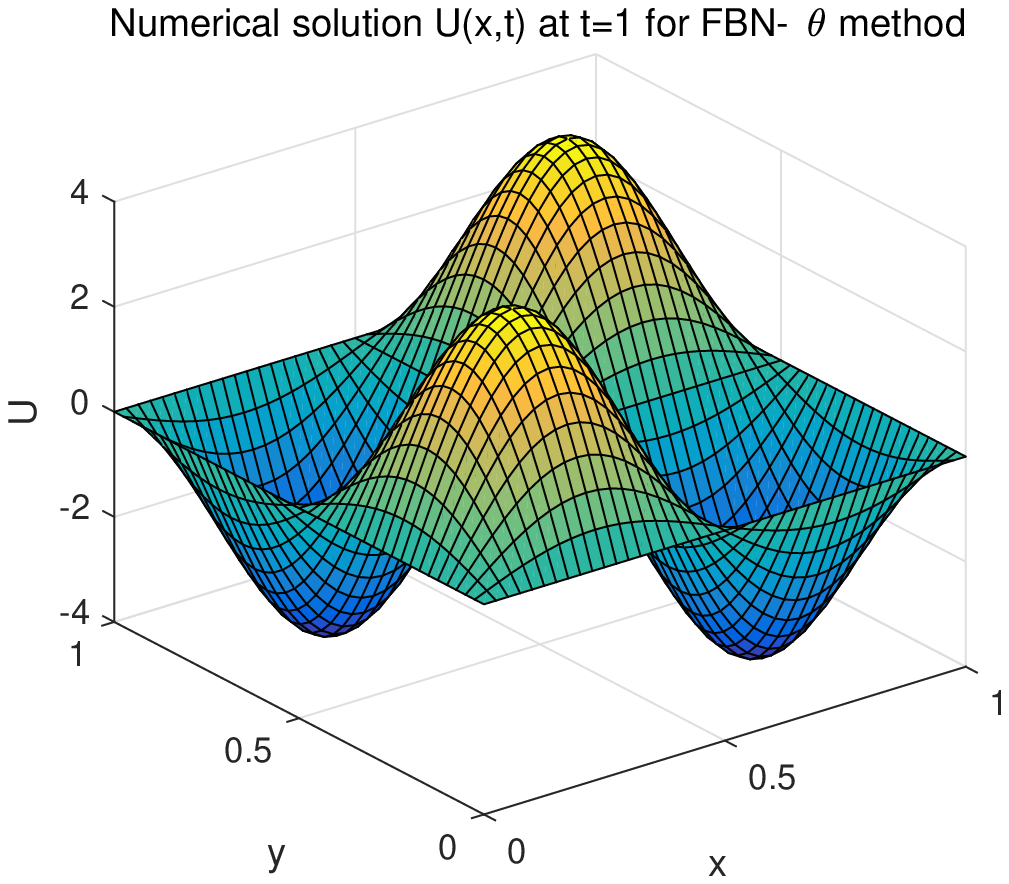}
  \caption{$h=\frac{\sqrt{2}}{40}$, $\tau=\frac{1}{200}$, $\gamma=0.4$, $\kappa=0.3$, $\theta_{\gamma}=-0.8$, $\theta_{\kappa}=1$.}\label{C6}
\end{minipage}
\end{center}
\end{figure}

In Table \ref{tab5} and Table \ref{tab6}, we test the convergence rate in time with fixed fine space mesh $h=\frac{\sqrt{2}}{400}$ for the FBT-$\theta$ and FBN-$\theta$ method, respectively.
Since the solution is sufficiently smooth for this example, we approximate the fractional derivative by the convolution part merely.
From the results of the tables, one can see that despite the different choices of fractional derivative orders $(\gamma,\kappa)$ and parameters $(\theta_{\gamma},\theta_{\kappa})$, the optimal second-order convergence rate is obtained.
To examine the spatial convergence rate, we fix the time step size $\tau=\frac{1}{200}$ in Table \ref{tab7} and Table \ref{tab8}, take $h=\frac{\sqrt{2}}{10},\frac{\sqrt{2}}{20},\frac{\sqrt{2}}{40}$, respectively.
The results of second-order convergence in space are in line with our theoretical analysis.
To intuitively compare the numerical solution with the exact one, we depict in Fig. \ref{C7} and Fig. \ref{C6} the numerical solutions at $t=T$ obtained by the FBT-$\theta$ and FBN-$\theta$ methods, respectively.
The space-time mesh is chosen as $h=\frac{\sqrt{2}}{40}, \tau=\frac{1}{200}$.
One can see that both of the FBT-$\theta$ and FBN-$\theta$ methods approximate the fractional derivative well by the comparison with the exact solution in Fig. \ref{C8}, which is based on the same space-time mesh.
\begin{table}[h]
\centering
 \caption{The temporal convergence rate for the FBT-$\theta$ method with $h=\frac{\sqrt{2}}{400}$}\label{tab5}
 \begin{tabular}{ccccc}
\toprule
  $(\gamma,\kappa)$ &$(\theta_{\gamma},\theta_{\kappa})$&  $\tau$ & $E_o(\tau, h)$ & rate \\
\midrule
	&	(0,0)	&	  1/10  	&	5.80147E-03	&	--	\\
	&		&	  1/20  	&	1.47696E-03	&	1.97 	\\
	&		&	  1/40  	&	3.47434E-04	&	2.09 	\\
(0.8,0.9)	&	(0,0.49)	&	  1/10  	&	5.77837E-03	&	--	\\
	&		&	  1/20  	&	1.47083E-03	&	1.97 	\\
	&		&	  1/40  	&	3.45859E-04	&	2.09 	\\
	&	(-0.5,0)	&	  1/10  	&	9.30248E-03	&	--	\\
	&		&	  1/20  	&	2.46264E-03	&	1.92 	\\
	&		&	  1/40  	&	6.07297E-04	&	2.02 	\\
\midrule
	&	(0.4,-0.1)	&	  1/10  	&	4.10911E-03	&	--	\\
	&		&	  1/20  	&	1.01088E-03	&	2.02 	\\
	&		&	  1/40  	&	2.26175E-04	&	2.16 	\\
(0.7,0.3)	&	(0.3,-1.5)	&	  1/10  	&	5.79435E-03	&	--	\\
	&		&	  1/20  	&	1.46303E-03	&	1.99 	\\
	&		&	  1/40  	&	3.42361E-04	&	2.10 	\\
	&	(-1,0)	&	  1/10  	&	1.85514E-02	&	--	\\
	&		&	  1/20  	&	5.10485E-03	&	1.86 	\\
	&		&	  1/40  	&	1.30957E-03	&	1.96 	\\
\bottomrule
 \end{tabular}
\end{table}
\begin{table}[h]
\centering
 \caption{The temporal convergence rate for the FBN-$\theta$ method with $h=\frac{\sqrt{2}}{400}$}\label{tab6}
 \begin{tabular}{ccccc}
\toprule
  $(\gamma,\kappa)$ &$(\theta_{\gamma},\theta_{\kappa})$&  $\tau$ & $E_o(\tau, h)$ & rate \\
\midrule
	&	(0,0)	&	  1/10  	&	2.21254E-02	&	--	\\
	&		&	  1/20  	&	5.72476E-03	&	1.95 	\\
	&		&	  1/40  	&	1.43163E-03	&	2.00 	\\
(0.2,0.8)	&	(0,0.5)	&	  1/10  	&	2.21130E-02	&	--	\\
	&		&	  1/20  	&	5.72128E-03	&	1.95 	\\
	&		&	  1/40  	&	1.43071E-03	&	2.00 	\\
	&	(0,1)	&	  1/10  	&	2.21403E-02	&	--	\\
	&		&	  1/20  	&	5.72857E-03	&	1.95 	\\
	&		&	  1/40  	&	1.43259E-03	&	2.00 	\\
\midrule
	&	(-1,-0.5)	&	  1/10  	&	5.91366E-02	&	--	\\
	&		&	  1/20  	&	1.61172E-02	&	1.88 	\\
	&		&	  1/40  	&	4.17263E-03	&	1.95 	\\
(0.5,0.6)	&	(-1,0.5)	&	  1/10  	&	5.89600E-02	&	--	\\
	&		&	  1/20  	&	1.60688E-02	&	1.88 	\\
	&		&	  1/40  	&	4.15997E-03	&	1.95 	\\
	&	(-1,1)	&	  1/10  	&	5.90674E-02	&	--	\\
	&		&	  1/20  	&	1.60965E-02	&	1.88 	\\
	&		&	  1/40  	&	4.16706E-03	&	1.95 	\\
\bottomrule
 \end{tabular}
\end{table}
\begin{table}[h]
\centering
 \caption{The spatial convergence rate for the FBT-$\theta$ method with $\tau=\frac{1}{200}$}\label{tab7}
 \begin{tabular}{ccccc}
\toprule
  $(\gamma,\kappa)$ &$(\theta_{\gamma},\theta_{\kappa})$&  $h$ & $E_o(\tau, h)$ & rate \\
\midrule
	&	(0,0)	&	  $\sqrt{2}$/10  	&	7.58676E-02	&	--	\\
	&		&	  $\sqrt{2}$/20  	&	1.89074E-02	&	2.01 	\\
	&		&	  $\sqrt{2}$/40  	&	4.71383E-03	&	2.00 	\\
(0.8,0.4)	&	(0.1,0.45)	&	  $\sqrt{2}$/10  	&	7.58695E-02	&	--	\\
	&		&	  $\sqrt{2}$/20  	&	1.89095E-02	&	2.00 	\\
	&		&	  $\sqrt{2}$/40  	&	4.71597E-03	&	2.00 	\\
	&	(-1,-2)	&	  $\sqrt{2}$/10  	&	7.58493E-02	&	--	\\
	&		&	  $\sqrt{2}$/20  	&	1.88881E-02	&	2.01 	\\
	&		&	  $\sqrt{2}$/40  	&	4.69437E-03	&	2.01 	\\
\bottomrule
 \end{tabular}
\end{table}
\begin{table}[h]
\centering
 \caption{The spatial convergence rate for the FBN-$\theta$ method with $\tau=\frac{1}{200}$}\label{tab8}
 \begin{tabular}{ccccc}
\toprule
  $(\gamma,\kappa)$ &$(\theta_{\gamma},\theta_{\kappa})$&  $h$ & $E_o(\tau, h)$ & rate \\
\midrule
	&	(0,0)	&	  $\sqrt{2}$/10  	&	7.54590E-02	&	--	\\
	&		&	  $\sqrt{2}$/20  	&	1.87869E-02	&	2.01 	\\
	&		&	  $\sqrt{2}$/40  	&	4.66558E-03	&	2.01 	\\
(0.4,0.3)	&	(0.5,0.5)	&	  $\sqrt{2}$/10  	&	7.54639E-02	&	--	\\
	&		&	  $\sqrt{2}$/20  	&	1.87923E-02	&	2.01 	\\
	&		&	  $\sqrt{2}$/40  	&	4.67102E-03	&	2.01 	\\
	&	(-0.8,1)	&	  $\sqrt{2}$/10  	&	7.53721E-02	&	--	\\
	&		&	  $\sqrt{2}$/20  	&	1.86908E-02	&	2.01 	\\
	&		&	  $\sqrt{2}$/40  	&	4.56767E-03	&	2.03 	\\
\bottomrule
 \end{tabular}
\end{table}
\section{Conclusion}
Two families of novel fractional $\theta$-methods are applied to approximate the fractional derivatives in the fractional Cable model.
With the help of the positivity properties of the coefficients of the methods, stability estimates and optimal convergence rate are derived.
For the case with solutions of weak regularity, the starting part is added to restore the second-order convergence rate in time.
\par
Nonetheless, the stability analysis and error estimates for the resulted scheme by the fractional $\theta$-methods are difficult for PDEs without the first derivative.
Authors think one reason is that to devise effective rules for the fractional $\theta$-methods with arbitrary parameter $\theta$ is extremely difficult.
A systematic approach for the analysis of the fractional $\theta$-methods when applied to PDEs without the first derivative is our future work.
\begin{acknowledgements}
The authors are grateful to Professor Buyang Li, two anonymous referees and editors for their valuable suggestions which improve the presentation of this work.
The work of the second author was supported in part by the NSFC grant 11661058.
The work of the third author was supported in part by the NSFC grant 11761053, the NSF of Inner Mongolia 2017MS0107, and the program for Young Talents of Science and Technology in Universities of Inner Mongolia Autonomous Region NJYT-17-A07.
The work of the fourth author was supported in part by grants NSFC 11871092 and U1930402.
\end{acknowledgements}

%
 \section*{Conflict of interest}

 The authors declare that they have no conflict of interest.



\end{document}